\input amstex
\documentstyle{amsppt}
\refstyle {A}
\magnification=\magstep1 

\hoffset .25 true in
\voffset .2 true in

\hsize=6.1 true in
\vsize=8.5 true in

\def\op{\operatorname}

\def\a{\alpha}
\def\b{\beta}

\def\d{\delta}
\def\l{\lambda}

\def\D{\Delta}
\def\HD{\hat\Delta}

\def\A{\bold A}

\def\F{\Bbb F}

\def\M{\Cal {M}_c}

\def\PP{\Cal P}
\def\RR{\Cal R}
\def\SS{\Cal S}

\def\S{\bold S}
\def\SH{\frak {Sh}_c}
\def\T{\bold T}

\def\op{\operatorname}
\def\lrar{\leftrightarrow}

\def\sub{\subset}
\def\subeq{\subseteq}

\def\Id{\operatorname{Id}}
\def\de{\operatorname{def}}

\def\hx{\hat X}

\def\hook{\hookrightarrow}

\topmatter
\title
Constructible sheaves on simplicial complexes\\
and Koszul duality 
\endtitle

\author
Maxim Vybornov\\
\endauthor

\address Department of Mathematics, Yale University, 
10 Hillhouse Ave,  New Haven, CT 06520 
\endaddress
\email  mv\@math.yale.edu
\endemail

\abstract We obtain a linear algebra data presentation of the 
category $\SH (X,\d)$ of constructible 
with respect to perverse triangulation sheaves on a finite 
simplicial complex 
$X$. We also establish Koszul duality between $\SH (X,\d)$ and the
category $\M (X,\d)$ of perverse sheaves constructible with 
respect to the triangulation
\endabstract

\endtopmatter

\document
\head Introduction \endhead

Let $X$ be a finite simplicial complex. There is a well known
linear algebra data description of (constructible 
with respect to the triangulation) sheaves of vector spaces
on $X$. A sheaf corresponds
to a gadget (called cellular sheaf) 
which assigns vector spaces to simplices
and linear maps to pairs of incident simplices.

Let  $\d$ be some perversity, that is 
a function $\d :{\Bbb Z}_{\geq 0}\to 
{\Bbb Z}$ satisfying some additional properties. 
R. MacPherson [M93] obtained a linear algebra data
descripton of the category of 
(constructible with respect to the triangulation) 
$\d$-perverse sheaves of vector spaces
on $X$, generalizing the above description. 
MacPherson's linear algebra gadgets are called 
cellular perverse sheaves. The category of cellular 
perverse sheaves is denoted by $\PP(X,\d)$. 

Given $X$ and a perversity $\d$ we can construct certain 
subsets of $X$ called perverse simplices. The first version 
of such subsets was introduced in [GM80]. 
One could say that the role of 
perverse simplices in the intersection homology 
theory is the same
as the role of usual simplices in the usual homology theory.
In this paper we introduce the category of constructible 
with respect to $(-\d)$-perverse simplices sheaves of vector 
spaces on $X$, which we call $\d$-sheaves. 
If $\d(k)=-k$, $k\geq 0$, is the bottom perversity,
$(-\d)$-perverse simplices are just usual simplices, so 
we obtain usual constructible sheaves.
We believe that 
$\d$-sheaves should be as natural for the study of intersection
homology as classical constructible 
sheaves are for the study of classical homology.
Our main result (Theorem B) gives a description of the category
of $\d$-sheaves in terms
of linear algebra data. Our linear 
algebra gadgets assign vector spaces to perverse simplices
and linear maps to pairs of ``incident'' perverse  simplices.
The category of such gadgets is denoted by $\RR (X,\d)$.

It is useful to identify a category of perverse sheaves 
or a category of linear algebra data with
the category of modules over some underlying algebra.
The algebra underlying the category $\PP(X,\d)$ (resp. 
$\RR (X,\d)$) is denoted by $A(X,\d)$ (resp. $B (X,\d)$). 
In many important cases underlying algebras turn out
to be Koszul. Koszulity was established for
algebras underlying perverse sheaves (middle perversity) on 
certain algebraic varieties in [BGSo] and [PS]. 
The Koszulity of $A(X,\d)$ and $B(X,\d)$ was conjectured
by the author and established by A. Polishchuk in [P].
Moreover, it was shown in [V] that the categories
$\PP(X,\d)$ and $\RR (X,\d)$ are (in some sense) 
Koszul dual to each other. This 
implies the Koszul duality between the category of perverse 
sheaves and the category of $\d$-sheaves. Therefore,
we could say that the two categories carry the same amount
of information about the topology of $X$.

We believe that the category of $\d$-sheaves should be considered
for more complex stratified objects and that it may prove to be
as important for the study of singular spaces as the
corresponding category of perverse sheaves. The principal
technical advantage of considering $\d$-sheaves vs. perverse 
sheaves is that we do not have to deal with the derived category 
explicitly. The price we have to pay for this is more complicated 
topology.
It may be said that replacing perverse sheaves by $\d$-sheaves
we move the intrinsic complexity of the theory from homological 
algebra to topology.

The paper is organized as follows. 
In part 1 we recall some notions from 
[GM80], [M93], [V] and [BBD].
In part 2 we state the theorems.
In part 3 we prove our main theorem.

Throughout the paper $\F$ stands for a (commutative) field,
$\op{char}\F=0$. All vector spaces are considered to be over 
$\F$.

\head 1. Preliminaries \endhead
\subhead  1.1.  Perversities and Perverse skeleta \endsubhead

Let $X$ be a connected finite simplicial complex.  Let $\hx$
be the first barisentric subdivision of $X$. If 
two simplices $\D$ and $\D'$ of $X$ are incident we write 
$\D\lrar\D'$.
We assume that
\roster
\item $\dim X=n$
\item $X$ is a regular cell complex
\item For every set of $0$-dimensional simplices 
$\{v_0,v_1\dots,v_r\}$
there is either exactly zero or exactly one $r$-dimensional
simplex $\D=\{v_0,v_1\dots,v_r\}$ with 
$\{v_0,v_1\dots,v_r\}$ as its vertices
\item all simplices are nondegenerate
\endroster
\definition{Definition 1.1.1([M93])} 
A {\sl perversity $\d :{\Bbb Z}_{\geq 0}\to 
{\Bbb Z}$} is a function from the
non-negative integers ${\Bbb Z}_{\geq 0}$ to the 
integers such that 
$\delta (0)=0$ and $\delta$ takes every interval \linebreak 
$\{0,1,\ldots ,k\}\subset {\Bbb Z}_{\geq 0}$ bijectively to 
an interval
$\{a,a+1,\ldots ,a+k\}\subset {\Bbb Z}$ for some 
$a\in {\Bbb Z}_{\leq 0}$.
In other words, a perversity is such a function $\delta$ that 
$\delta(0)=0$ and for $k\in\Bbb Z_{\geq 0}$,
$$\delta(k)=
\cases \text {either} & \underset{i\in[0,k-1]}
\to{\max} \delta(i)+1 \\ 
\text { or} & \underset{i\in[0,k-1]}\to{\min} 
\delta(i)-1 \\ \endcases 
$$
\enddefinition
\definition{Definition 1.1.2([M93], cf. [GM80])} 
Given a perversity $\d$, 
we define $\d(\D)=\d(\dim\D)$, where $\D$ is a simplex of $X$.
Given a simplex $\HD=\{c_1,c_2,\dots,c_s\}$ of $\hx$ where 
$\{c_1,c_2,\dots,c_s\}$ are baricenters of $\D_1, \D_2,\dots,\D_s$,
$\dim\D_1<\dim\D_2<\dots<\dim\D_s$ we denote by $\max\D$ such
vertex $c_i$ that $\d(\D_i)=
\max\{\d(\D_1),\d(\D_2),\dots, \d(\D_s)\}$.
Given a simplex $\D$ with the baricenter $c$
we define the corresponding {\sl perverse simplex}
$$
^\d\!\D=\bigsqcup_{\max\HD=c} \HD
$$
We define the {\sl $k$-th perverse skeleton} $X^{\d}_k$, 
$\min_{[0,n]}\d\leq k\leq\max_{[0,n]}\d$, as follows
$$
X^{\d}_k=\bigsqcup_{\d(\D)\leq k}{}^\d\!\D\sub X
$$
Thus, we have a filtration
$$
X^{\d}_i\sub X^{\d}_{i+1}\sub\dots\sub X^{\d}_{i+n-1}
\sub X^{\d}_{i+n}=X
$$
where $i=\min_{[0,n]}\d$.
It is easy to see that perverse simplices are connected 
components of $X^{\d}_k-X^{\d}_{k-1}$, 
$\min_{[0,n]}\d\leq k\leq\max_{[0,n]}\d$. The decomposition of 
$X$ into a disjoint union of perverse simplices is called
{\sl $\d$-perverse triangulation} of $X$.
\enddefinition
From now on we will fix a perversity $\d$.
\subhead  1.2. Quiver algebras and categories of their modules 
\endsubhead

In this section we recall some definition from [V]. 
\definition{Definition 1.2.1} 
Let $\op{Q}(X,\d)$ be a quiver (i.e. finite simple oriented tree) 
whose vertices
are indexed by simplices of $X$ and there is an arrow from $\D$
to $\D'$ if and only if $\d(\D)=\d(\D')+1$ and $\D\lrar\D'$. 
There is a standard construction of a {\sl quiver algebra}
$\F\op{Q}(X,\d)$ associated to $\op{Q}(X,\d)$.
\enddefinition
\definition{Definition 1.2.2} 
Let $A(X,\d)$ be the quotient of 
$\F \op{Q}(X,\d)$ by {\bf the chain complex relations:}
if $\D'$ is any simplex such that
$\d(\D')=k+1$ and $\D''$ is any simplex such that 
$\d(\D'')=k-1$, then
$$\sum\limits_{\matrix
\D : \d(\D)=k,\\
\D'\lrar\D\lrar\D''\endmatrix}
a(\D ,\D'')\cdot a(\D',\D )=0$$
where $a(\D_\a, \D_\b)$ are generators of $A(X,\d)$.
\enddefinition

\definition{Definition 1.2.3}
Let $B(X,\d)$  be the quotient of
$\F \op{Q}(X,\d)$ by {\bf the equivalence relations:}
Suppose that $\D'$, $\D''$, $\D_1$ and $\D_2$ are cells
of $X$  such that:
\roster
\item $\d(\D')=k+1$, $\d (\D'')=k-1$ and $\d(\D_1)=\d(\D_2)=k$ 
\item  $\D'\lrar\D_1\lrar\D''$ and $\D'\lrar\D_2\lrar\D''$
\endroster
then 
$$
b(\D_1, \D'')\cdot b(\D',\D_1)=
b(\D_2, \D'')\cdot b(\D',\D_2)
$$
where $b(\D_\a, \D_\b)$ are generators of $B(X,\d)$.
\enddefinition 
The category of left finite dimensional modules over $A(X,\d)$
will be denoted by $\PP (X,\d)$ and 
the category of left finite dimensional modules over $B(X,\d)$
will be denoted by $\RR (X,\d)$. The category $\PP (X,\d)$
was introduced by R. MacPherson in [M93] for an
arbitrary finite regular cell complex. It is called the category
of {\sl cellular perverse sheaves}.
\subhead  1.3. Constructible sheaves \endsubhead
\definition{Definition 1.3.1} A sheaf $\A$ of $\F$-vector spaces
is called {\sl constructible}
with respect to $\d$-perverse triangulation if
for all simplices $\D$, 
$i_\D^*\A$ is a constant sheaf on ${}^\d\!\D$ 
associated to a finite dimensional vector space over $\F$, 
where $i_\D:{}^\d\!\D\hook X$
is the inclusion of the corresponding perverse simplex $^\d\!\D$
\enddefinition
We will denote the category of all constructible
with respect to $(-\d)$-perverse triangulation sheaves 
($\d$-sheaves) by $\SH (X,\d)$. 

\definition{Definition 1.3.2([GM80])} 
A {\sl "classical" perversity  $p :{\Bbb Z}_{\geq 0}\to 
{\Bbb Z}_{\geq 0}$} is a non-decreasing function 
from the set of
non-negative integers ${\Bbb Z}_{\geq 0}$ to itself such that 
$p (0)=0$ and 
$p(k)-p(k-1)$ is either $0$ or $1$. 
\enddefinition

It is easy to see that there is a one-to-one correspondence
between the set of all perversities and the set of all
"classical" perversities. Let $p$ be a classical perversity
corresponding to our fixed perversity $\d$. We will
denote by $\M (X,\d)$ the category of 
(homologically) constructible with respect to the triangulation
{\sl $p$-perverse sheaves} introduced in [BBD]. 
\head 2. Theorems \endhead

\proclaim{Theorem A} The following categories
$$
\RR (X)\simeq\SH (X)
$$
are (canonically) equivalent
\endproclaim
This result is rather well known.
A version of it was presented in [M93], 
another version is offered as exercise 
VIII.1 in [KS]. A detailed discussion will appear in [GMMV].
The primary goal of this paper is to generalize this result 
to the case of an arbitrary perversity
\proclaim{Theorem B} The following categories
$$
\RR (X,\d)\simeq\SH (X,\d)
$$
are (canonically) equivalent
\endproclaim
{\sl Canonically} equivalent means that the equivalence 
functor does 
not depend on anything but the triangulation of $X$ and 
perversity.
\demo{Proof} The proof will appear in part 3.\qed\enddemo

The following result was obtained by R. MacPherson ([M93], [M95]).
It may be considered as another generalization  of the Theorem A. 
 
\proclaim{Theorem C} The following categories 
$$
\PP (X,\d)\simeq\M (X, \d)
$$
are equivalent.
\endproclaim
\demo{Proof} It follows from [P] that $\M (X, \d)$ is equivalent
to the category of modules over the quadratic dual
to $B(X,\d)^{\text{opp}}$. Such a quadratic dual algebra is 
exactly $A(X,\d)$
\qed\enddemo

\proclaim{Theorem D (cf. [P], [V])} 
(a) The category $\SH (X,\d)$ is Koszul;

(b) The category $\M (X,\d)$ is Koszul;

(c) There exists a functor
$$K:D^b(\M (X,\d))\to D^b(\SH (X,-\d)) $$

which is an equivalence of triangulated categories.

\endproclaim
\demo{Proof} The proof of the fact that  $B(X,\d)$ 
is Koszul (up to changing the right multiplication
to the left one) is given in 
section 3 of [P]. Then 
(a) is implied by Theorem B.
Since the quadratic dual $A(X,\d)^!=B(X, -\d)$ it follows
that $A(X,\d)$ is Koszul as the quadratic dual of a Koszul algebra.
Therefore (b) is implied by Theorem C. Quadratic duality
of $A(X,\d)$ and $B(X, -\d)$ and the fact that 
$B(X,\d)=B(X, -\d)^{\text{opp}}$ implies that $E(A(X,\d))=B(X,\d)$
by Theorem 2.10.1 of [BGSo], where $E(A(X,\d))$ is the Koszul
dual of $A(X,\d)$. Using the construction of [BGSo] it is not hard 
to see (cf. [V]) that the Koszul duality 
functor restricts to an equivalence functor 
$K: D^b(\PP (X,\d))\to D^b(\RR (X,-\d))$.
Thus, (c) follows from Theorem B and Theorem C. 
\qed\enddemo
\remark{Remark} Note that the Koszul duality functor 
transforms simple objects in $\M (X,\d)$ (intersection 
homology sheaves on closed simplices) to
indecomposable projective objects in $\SH (X,\d)$ 
(constant sheaves on certain subsets of $X$).
\endremark

\head 3. Proof of the Theorem B \endhead 
\subhead 3.1. Linear algebra data \endsubhead

We will redefine the category $\RR (X, \d)$ in terms of linear 
algebra data. We leave it to the reader to check that the two
definitions coincide.
\definition{Definition 3.1.1}
An object $\bold S$ of $\RR (X,\d)$
is the following data:
\roster
\item  (Stalks) For every simplex $\D$ in $X$, a finite dimensional
vector space 
$S(\D)$ called the {\sl stalk} of ${\S}$ at $\D$
\item  (Restriction  maps) For every pair of simplices $\D$ and 
$\D'$ in 
$X$ such that $\d(\D)=\d(\D')+1$, and $\D\lrar\D'$, a linear map
$s(\D ,\D'):S(\D)\to S(\D')$ called the 
{\sl restriction map}
\endroster
subject to {\bf the equivalence axiom (TEA)}:  

Suppose that $\D'$, $\D''$, $\D_1$ and $\D_2$ are simplices
of $X$  such that:
\roster
\item $\d(\D')=k+1$, $\d (\D'')=k-1$ and $\d(\D_1)=\d(\D_2)=k$ 
\item $\D'\lrar\D_1\lrar\D''$ and $\D'\lrar\D_2\lrar\D''$
\endroster
then
$$
s(\D_1, \D'')\circ s(\D',\D_1)=
s(\D_2, \D'')\circ s(\D',\D_2)
$$
The morphisms in this category are stalkwise linear maps 
commuting with
the restriction maps.
\enddefinition

We define another category $\SS (X,\d)$ of linear algebra data. 
\definition{Definition 3.1.2} 
The category $\SS (X,\d)$
is a full subcategory of the abelian category $\RR (\hx)$. 
An object $\T$ of $\RR (\hx)$
belongs to $\SS(X,\d)$ 
if for any $\HD,\HD'\subseteq {^{-\d}\!\D}$ we have
\roster
\item $T(\HD)=T(\HD')$
\item if $\HD$ is a codim-$1$ face of $\HD'$ then
$t(\HD,\HD')=\op{Id}_{T(\HD)}$.
\endroster
\enddefinition

\subhead 3.2. $\RR (X,\d)=\SS (X,\d)$ \endsubhead

\definition{Definition 3.2.1} We will consider a partial order on 
the set of simplices of $X$ which is uniquely defined by the
following: $\D<\D'$ if 
\roster
\item $\D$ and $\D'$ are incident
\item $\d(\D)=\d(\D')+1$
\endroster
we say that $\D\leq\D'$ if $\D<\D'$ or $\D=\D'$. It is easy to see
that if $\D\lrar\D'$ and $\d(\D)>\d(\D')$ then $\D<\D'$.
\enddefinition

\definition{Definition 3.2.2} Let $\D_I\leq\D_T$ be two simplices
of $X$. Let $\S$ be an object of $\RR (X,\d)$. We define a linear map
$f(\D_I, \D_T): S(\D_I)\to S(\D_T)$ in the following way
\roster
\item if  $\D_I=\D_T$ then $f(\D_I, \D_T)=\op{Id}_{S(\D_I)}$
\item if $\D_I<\D_T$ then by definition there is a sequence
$$\D_I=\D_0<\D_1<\D_2<\dots<\D_r=\D_T\tag 3.2.1$$
$\d(\D_i)=\d(\D_{i+1})+1$, $0\leq i\leq r-1$. We set
$$f(\D_I, \D_T)=s(\D_{r-1},\D_r)\circ\dots\circ s(\D_0,\D_1)$$
\endroster
It is easy to see that $f(\D_I, \D_T)$ does not depend on the
choice of sequence (3.2.1); and moreover, that all $<$ in 
(3.2.1) could be changed to $\leq$. In particular, 
if $\D_1\leq\D_2\leq\D_3$ then
$$
f(\D_1,\D_3)=f(\D_2,\D_3)\circ f(\D_1,\D_2)\tag 3.2.2
$$ 
\enddefinition

\proclaim{Lemma 3.2.3} 
Let $\hat\D,\hat\D'$ of $\hat X$
be such that  $\hat\D$ is a $\op{codim}\ 1$ face of $\hat\D'$.
Let $\HD\subseteq{^{-\d}\!\D}$,
$\hat\D'\subseteq{^{-\d}\!\D'}$. Then
$\D\lrar\D'$ and $\D\leq\D'$.
\endproclaim
\demo{Proof} Let
$\hat\D=\{v_1,\dots,v_k\}$, $\hat\D'=\{v_1,\dots,v_k,e\}$.
Case I: $e=\max(\HD')$.  Then $e$ is a baricenter of $\D'$.
By definitions $-\d(\D')>-\d(\D)\implies\d(\D')<\d(\D)$.
Let $v_i=\max(\HD)$, $v_i$ is the baricenter of $\D$. Then
1-simplex $\{v_i, e\}$ is a simplex of $\hx$, thus $\D\lrar\D'$.
Case II: $e\neq \max(\hat\D')$. Then $v_i=\max(\HD)=\max(\HD')$,
$v_i$ is the baricenter of both $\D$ and $\D'$, thus $\D=\D'$.
\qed\enddemo

\proclaim{Theorem 3.2.4} The category $\RR (X,\d)$ is 
isomorphic
to the category $\SS (X,\d)$
$$
\RR (X,\d)=\SS (X,\d)
$$
\endproclaim
\demo{Proof} (a) The functor $\Phi:\RR (X,\d)\to \SS(X,\d)$.
If $\S$ is an object of $\RR (X,\d)$ then
$\T=\Phi(\S)$ is constructed as follows

$$
\matrix\format\l & \l \\
T(\HD)=S(\D), & \text {   for } \HD\subseteq{}^{-\d}\!\D \\
\endmatrix
$$
If $\hat\D'$ is a $\op{codim}\ 1$ face of $\hat\D''$ then
$\D'\leq\D''$ by Lemma 3.2.3 and we set
$$
\matrix\format\l & \l \\
t(\HD', \HD'')=f(\D', \D''), & \text {   for } 
\HD'\subseteq{}^{-\d}\!\D' \text { and } \HD''\subseteq{}^{-\d}\!\D''
\endmatrix
$$
Now let $\HD'\subeq{}^{-\d}\!\D'$, $\HD_1\subeq{}^{-\d}\!\D_1$,
$\HD_2\subeq{}^{-\d}\!\D_2$, $\HD''\subeq{}^{-\d}\!\D''$ be such
a quadruple for which we have to check TEA. 
Then by Lemma 3.2.3 we have
$$
\matrix
\D'\leq\D_1\leq\D'' \\
\D'\leq\D_2\leq\D'' \\
\endmatrix
$$
By definitions  we have
$$
\aligned
t(\HD_1, \HD'')\circ t(\HD', \HD_1)
\overset{\de}\to= & f(\D_1,\D'')\circ f(\D',\D_1) 
\overset{(3.2.2)}\to= f(\D',\D'')\\
\overset{(3.2.2)}\to= &f(\D_2,\D'')\circ f(\D',\D_2)
\overset{\de}\to=t(\HD_2, \HD'')\circ t(\HD', \HD_2)
\endaligned
$$
$\Phi$ on morphisms is defined in the obvious way. 

(b) The functor $\SS (X,\d)@>>>\RR (X,\d)$.
If $\T$ is an object of $\SS (X,\d)$ then
$\S=\Psi(\T)$ is constructed as follows

$$
S(\D)=T(\HD),\qquad \text {for } \HD\subseteq {}^{-\d}\!\D 
$$
$S(\D)$ is well defined due to constructibility of $\T$.
Let $\D'$ and $\D''$ be two incident simplices of $X$
such that $\d(\D')=\d(\D'')+1$. We have to construct
the map $s(\D',\D'')$. Let $c'$ be a baricenter of $\D'$
and $c''$ be a baricenter of $\D''$. We set
$$
s(\D',\D'')=t(c',\{c',c''\})
$$
where $\{c',c''\}$ is a simplex of $\hx$.
Let $\D'<\D_m<\D''$, $\d(\D')=\d(\D_m)+1=\d(\D'')+2$.
Let $c'$, $c_m$ and $c''$ be baricenters of $\D'$, $\D_1$
and $\D''$ respectively. Let us assign special names to
the following four simplices of $\hx$: 
$l=\{c',c''\}$, $l_m'=\{c',c_m\}$, $l_m''=\{c'',c_m\}$
and $tr=\{c',c_m, c''\}$. We have
$$
\aligned
s(\D',\D'')\overset{\de}\to= & t(c',l) 
\overset{t(l, tr)=\Id}\to = t(l, tr)\circ t(c',l) 
\overset{\text{TEA}}\to= 
t(l_m', tr)\circ t(c', l_m') \\
\overset{t(c_m,l_m')=\Id}\to = &
t(l_m', tr)\circ t(c_m,l_m')\circ t(c',l_m')  
\overset{\text{TEA}}\to= 
t(l_m'', tr)\circ t(c_m,l_m'')\circ t(c',l_m')\\ 
\overset{t(l_m'', tr)=\Id }\to = &
t(c_m,l_m'')\circ t(c',l_m')
\overset{\de}\to = s(\D_m,\D'')\circ s(\D',\D_m)
\endaligned
$$
The equivalence axiom follows.

$\Psi$ on morphisms is defined in the obvious way.

(c) It follows from our explicit construction that
$\Psi\circ\Phi=\Id$. 
Using constructibility of $\T$, definitions and (3.2.2)
it is easy to see that
$\Phi\circ\Psi=\Id$.
\qed\enddemo

\subhead 3.3. $\SS (X,\d)\simeq\SH (X,\d)$ \endsubhead

\definition{Definition 3.3.1} Let $\S$ be an object of $\RR (\hx)$.
Let $i_A:A\hook X$ be a closed union of simplices of $\hx$. 
We define
a functor $i_A^*:\RR (\hx) \to \RR (A)$ as
follows: if $\T=i_A^*\S$ then 
\roster
\item $\T(\HD)=\S(\HD)$
\item $t(\HD', \HD'')=s(\HD', \HD'')$
\endroster
\enddefinition

\proclaim{Lemma 3.3.2} The following functorial diagram commutes
i.e. two possible compositions of functors are isomorphic
$$
\CD
\RR (\hx) @<\sim<< \SH (\hx) \\
@V i^*_A VV @V i^*_A VV \\
\RR (A)  @<\sim<<  \SH (A) \\
\endCD
$$
where $i^*_A: \SH (\hx)\to \SH (A)$ is a standard sheaf 
theory functor.
\endproclaim
\demo{Proof} Proof is left to the reader \qed\enddemo

\proclaim{Theorem 3.3.3} The category $\SS (X,\d)$ is 
(canonically) equivalent
to the category $\SH (X,\d)$
$$
\SS (X,\d)\simeq\SH (X,\d)
$$
\endproclaim
\demo{Proof} By definitions and functoriality of $i^*$ the
category $\SH (X,\d)$ is a full subcategory of $\SH (\hx)$.
Lemma 3.3.2 implies that the equivalence functor 
$\RR (\hx)\simeq\SH (\hx)$ restricts to an equivalence functor
$\SS (X,\d)\simeq\SH (X,\d)$.\qed\enddemo
Theorem 3.2.4 and Theorem 3.3.3 imply Theorem B.
\head Acknowlegdments \endhead

The primary motivation of my work stems from
questions formulated by M. Goresky and R. MacPherson.
I am very grateful to them for 
useful conversations and steady encouragement. I would also like
to acknowlegde my special indebtedness to Robert MacPherson, who
invented cellular perverse sheaves and from whom I learned
about them. 

\enddocument
\end